\documentclass[journal,9pt]{IEEEtran}
\ifCLASSINFOpdf
\else
\fi
\usepackage{epsfig}
\usepackage{amssymb}
\usepackage{amsmath}
\usepackage{amsfonts}
\usepackage{color}
\usepackage{epstopdf}
\usepackage{pdfsync}
\usepackage{multirow}
\usepackage{caption2}
\newtheorem{theorem}{Theorem}

\newtheorem{lemma}{Lemma}

\newtheorem{proof of Theorem 1}{Proof of Theorem 1}

\begin{document}
\title{ Adaptive Output Consensus of Heterogeneous Nonlinear Multi-agent Systems: A Distributed Dynamic Compensator Approach}
\author{Guangqi Li,  Long Wang
\thanks{This work was supported in part by National Natural Science Foundation of China under Grants 61751301 and 61533001, in part by the National Natural Science Foundation of Shaanxi Province under Grant 2019JQ-190, and in part by Fundamental Research Funds for the Central Universities under grant JB180406.}
 \thanks{Guangqi Li is with the School of Mechano-Electronic Engineering, Xidian University, Xian, China, 710071.
Long Wang is with the Center for Systems and Control, College of Engineering, Peking University, Beijing, China, 100871. }
\thanks{%
Corresponding author: Long Wang. E-mail addresses: gqli@xidian.edu.cn,  longwang@pku.edu.cn.
}}
\maketitle

\begin{abstract}
Distributed dynamic compensators, also known as distributed observer, play a key role in the output consensus problem of heterogeneous nonlinear multi-agent systems. However, most existing distributed dynamic compensators require either the compensators' information to be exchanged through communication networks, or that the controller for each subsystem satisfies a class of small gain conditions. In this note, we develop a novel distributed dynamic compensator to address the adaptive output consensus problem of heterogeneous nonlinear multi-agent systems with unknown parameters. The distributed dynamic compensator only requires the output information to be exchanged through communication networks. Thus, it reduces the communication burden and facilitates the implementation of the dynamic compensator. In addition, the distributed dynamic compensator converts the original adaptive output consensus problem into the global asymptotic tracking problem for a class of nonlinear systems with unknown parameters.
 Then, by using the adaptive backstepping approach, we develop an adaptive tracking controller for each subsystem, which does not requre  the small gain conditions as in previous studies. It is further proved that all signals in the
closed-loop system are globally uniformly bounded, and the proposed
scheme enables the outputs of all the subsystems to track the output
of leader asymptotically. A simulation is presented to illustrate the effectiveness of the design methodology.
\end{abstract}
\markboth{PREPRINT SUBMITTED TO IEEE TRANSACTIONS ON AUTOMATIC
CONTROL}{Shell \MakeLowercase{\textit{et al.}}: Bare Demo of
IEEEtran.cls for Journals}

\IEEEpeerreviewmaketitle
\section{Introduction}
The consensus problem of multi-agent systems has attracted many researchers, due to its widespread potential applications in various fields. Its objective is to design a distributed control law  such that the states or the outputs of all agents achieve an agreement. The control law is distributed in the sense that each agent's controller only uses information from the agent and its neighboring agents. During the past decades, the consensus problem for multi-agent systems has been extensively studied from various perspectives \cite{Jadbabaie-AC-2003}-\cite{ShiGuodong-SIAM-2013}.  For more details, please refer to the surveys \cite{Olfati-IEEE-2007}-\cite{Qin-IE-2017} and the references cited therein.\par
Recently, more attention has been paid to the heterogeneous nonlinear multi-agent systems \cite{Yang-IJACSP-2018}-\cite{Zhu-AC-2016}. Distributed dynamic compensators, also called distributed observer, are useful in dealing with the output consensus problem of heterogeneous nonlinear multi-agent systems. This problem can be addressed in two step. First, a local dynamic compensator is designed for each agent, and the outputs or states of all compensators achieve consensus through a proper collaborative control strategy. Then, the output regulation theory is applied to constructing controller, forcing the output of each agent to track the output of local compensator. Based on this method, the output consensus problem has been addressed for different classes of heterogeneous nonlinear multi-agent systems \cite{Isidori-AC-2014}-\cite{LiuTao-AC-2019}. For example, the output synchronization problem was investigated in \cite{Isidori-AC-2014} for heterogeneous nonlinear multi-agent systems.  The cooperative output regulation problem was addressed for heterogeneous nonlinear multi-agent systems with unknown and non-identical control directions \cite{Guo-AC-2017}. Unfortunately, all the distributed dynamic compensators in \cite{Isidori-AC-2014}-\cite{LiuTao-AC-2019} require the compensator information to be exchanged through communication networks. The compensator information is not physical but artificial, hence exchanging such information must incur additional communication complexity and burden. In many physical circumstances, each agent can only observe or measure the output information of its neighboring agents. As a result, it is more desirable to design distributed controller under output communication. However, the output communication also brings new challenges in designing controller, and new design technique is required. In a recent paper \cite{Zhu-AC-2016}, a general framework was proposed to address the output consensus problem of heterogeneous nonlinear multi-agent systems under output communication. Actually, the distributed dynamic compensator constructed under output communication and each subsystem can be viewed as a interconnection system. Then, the controller satisfying a class of small gain conditions is designed for each subsystem to address the tracking problem of the interconnection system. However, the small gain conditions result in sufficiently large control gains, and for some nonlinear systems with completely unknown parameters, it is unable or difficult to design controller satisfying small gain conditions.         \par
In this note, a novel distributed dynamic compensator is developed to address the adaptive output consensus problem for heterogeneous nonlinear multi-agent systems with unknown parameters. The distributed dynamic compensator only requires the output information to be exchanged through communication networks. This considerably reduces the communication burden and facilitates the implementation of the dynamic compensator. In addition, the distributed dynamic compensator converts the adaptive output consensus problem of heterogeneous nonlinear multi-agent systems with unknown parameters into the problem of global asymptotic tracking for a class of nonlinear systems with unknown parameters. Then, by using adaptive backstepping approach, we develop an adaptive tracking controller for each subsystem, without requiring the small gain conditions \cite{Zhu-AC-2016}. It is further proved that all signals in the closed-loop system are globally uniformly bounded, and the proposed scheme enables the outputs of all the subsystems to track the output of leader asymptotically. A simulation is presented to illustrate the effectiveness of the design methodology.\par
The adaptive output consensus problem has also been addressed via adaptive backstepping approach \cite{Kristic-NACD-1995} for nonlinear multi-agent systems with unknown parameters \cite{Wangwei-A-2014}-\cite{HuangJingshuai-A-2017}.
Compared with these results \cite{Wangwei-A-2014}-\cite{HuangJingshuai-A-2017}, our designed methodology has the following advantages:
\begin{itemize}
\item  In \cite{Wangwei-A-2014}-\cite{HuangJingshuai-A-2017}, some restrictive conditions were imposed, e.g., each agent needs to know the state information and nonlinear functions of its neighbors \cite{Wangwei-A-2014}-\cite{Wangwei-A-2017}, or the filter information of its neighbors \cite{HuangJingshuai-A-2017}, and the system orders of all agents need to be the same. However, in our work,  the system orders
are not the same for all agents, and only the output information is exchanged through communication netowrk.
\item Our design methodology is more flexible than the methods in \cite{Wangwei-A-2014}-\cite{HuangJingshuai-A-2017}. Actually, by means of the distributed dynamic compensator, we can use different control approaches to design tracking controller for each subsystem. Thus, our proposed methodology can be used to address the output consensus problems of heterogeneous nonlinear multi-agent systems with non-identical structure, and hence the output consensus problem of multi-agent systems with unknown and non-identical control directions. However, it is difficult to apply the methods in \cite{Wangwei-A-2014}-\cite{HuangJingshuai-A-2017} to these problems, even for the case that the system orders of all agents are the same.
\item In \cite{Wangwei-A-2014}-\cite{Wangwei-A-2017}, each agent required constructing additional local estimates to account for the unknown parameters of its neighbors' dynamics. This inevitably results in a much complex controller. However, in our work, each agent does not need to construct the additional local estimates.
\end{itemize}

The rest of this note is organized as follows. In Section 2, we formulate our problem statement, and give some useful lemmas. In Section 3, we first develop a novel distributed dynamic compensator to address the challenges caused by heterogeneous dynamics. Then, adaptive controller design and stability analysis of closed-loop system are presented.  Finally, an illustrative example is  provided in Section 4.\par

{\bf \textit{Notation.}} Throughout this note, $\mathbb{R}$ and $\mathbb{R}^m$ denote the set of real numbers and $m$-dimensional real vector space, respectively. $0_{r}$ denotes the column vector with all $r$ elements being $0$, and $I_{r\times r}$ denotes the $r\times r$ identity matrix. For a given vector function $z(t)\in\mathbb{R}^{m}$, $\|z(t)\|$ and $\|z(t)\|_{\infty}$  denote the standard Euclidean norm and the essential supremum norm, respectively. Moreover, a function $\gamma:\mathbb{R}^+\rightarrow\mathbb{R}^+$ is said to be a $\mathcal{K}$-class function if it is continuous, strictly increasing, and $\gamma(0)=0$; a function $\beta:\mathbb{R}^+\times \mathbb{R}^+\rightarrow \mathbb{R}$ is said to be a $\mathcal{KL}$-class function if $\beta(\cdot,t)$ is of class $\mathcal{K}$ for each fixed $t>0$ and $\beta(s,t)\rightarrow 0$ decreasingly as $t\rightarrow \infty$ for each fixed $s>0$.

\section{Problem formulation and Preliminaries}
In this note, we consider the leader-following output consensus problem for the following class of heterogeneous nonlinear multi-agent systems with unknown parameters:
\begin{eqnarray}
&&\dot{x}_{i,l}=x_{i,l+1}+\psi_{i,l}(x_{i,1},\cdots,x_{i,l})^T\theta_{i},l=1,\cdots,r_i-1,\nonumber\\
&&\dot{x}_{i,r_i}=u_i+\psi_{i,r_i}(x_{i,1},\cdots,x_{i,r_i})^T\theta_{i},\nonumber\\
&&y_i=x_{i,1},i=1,\cdots,N,\label{systemF}
\end{eqnarray}
where $x_i=[x_{i,1},\cdots,x_{i,r_i}]^T\in\mathbb{R}^{r_i},y_i\in\mathbb{R},u_i\in\mathbb{R}$ are the state, output and input of the $i$th subsystem, respectively, $\theta_i\in\mathbb{R}^{m_i}$ a vector of unknown constants, and $\psi_{i,l}(\cdot):\mathbb{R}^{l}\rightarrow \mathbb{R}^{m_i}$ known smooth nonlinear functions. The system is heterogeneous in the sense that the order $r_i$, and the nonlinear functions $\psi_{i,l}(\cdot)$ need not to be identical for all agents.\par
The leader's signal $y_0\in \mathbb{R}$ is assumed to be generated by
a linear autonomous system of the form
\begin{eqnarray}
&&\dot{x}_0=Ax_0,\nonumber\\
&&y_0=Cx_0,\label{systemL}
\end{eqnarray}
where $x_0\in\mathbb{R}^\nu$, $A\in\mathbb{R}^{\nu\times \nu}$ and $C\in\mathbb{R}^{1\times \nu}$. Without loss of generality, we assume that $(A,C)$ is detectable.\par
The communication network of this multi-agent system
can be described by a digraph $\bar{\mathcal{G}}=\{\bar{\mathcal{V}},\bar{\mathcal{E}}\}$ with  $\bar{\mathcal{V}}=\{v_0,v_1,\cdots,v_N\}$ and $\bar{\mathcal{E}}\subseteq \bar{\mathcal{V}}\times \bar{\mathcal{V}}$, where node $v_0$ is associated with the leader system $(\ref{systemL})$, and node $v_i,i=1,\cdots,N$ is associated with the $i$th subsystem of $(\ref{systemF})$. For $i=1,\cdots, N,j=0,\cdots, N$ and $i\neq j$, $(v_j,v_i)\in \bar{\mathcal{E}}$ if and only if the control law $u_i$ can use the output information of $j$th subsystem or leader for control. Let $\bar{\mathcal{A}}=[a_{ij}]\in\mathbb{R}^{(N+1)\times (N+1)}$ be the weighted adjacency matrix of $\bar{\mathcal{G}}$, where $a_{ii}=0,i=0,1,\cdots,N$, and $a_{ij}>0,i=1,\cdots,N,j=0,1,\cdots,N$ if and only if $(v_j,v_i)\in \bar{\mathcal{E}}$. The neighbor set of
agent $i$ is defined as $\bar{\mathcal{N}}_i=\{v_j,(v_j,v_i)\in\bar{\mathcal{E}}\}$. Let $\mathcal{G}=\{\mathcal{V},\mathcal{E}\}$ be the subgraph of $\bar{\mathcal{G}}$, where $\mathcal{V}=\{v_1,\cdots,v_N\}$, and $\mathcal{E}\subseteq \mathcal{V}\times \mathcal{V}$ is obtained from $\bar{\mathcal{E}}$ by removing all edges between the node $v_0$ and nodes in $\mathcal{V}$. \par
Let us describe our control law as follows:
\begin{eqnarray}
u_i&=&h_i(x_i,\xi_i),\nonumber\\
\dot{\xi}_i&=&l_i(\xi_i,x_i,e_{vi}),i=1,\cdots,N,\label{control}
\end{eqnarray}
where $h_i(\cdot)$ and $l_i(\cdot)$ are some nonlinear functions, and $e_{vi}=\sum_{j=0}^Na_{ij}(y_i-y_j)$.\par
A control law of the form $(\ref{control})$ is distributed  since $u_i$ only depends on the output information of its neighbor and the state information of itself. Our problem is described as follows.\\
\textbf{Problem 1.}~Given the multi-agent systems  $(\ref{systemF})$-$(\ref{systemL})$, and a digraph $\bar{\mathcal{G}}$, design a control law of the form $(\ref{control})$, such that the solution of the closed-loop system is globally uniformly bounded, and satisfies $\lim_{t\rightarrow \infty}(e_i:=y_i-y_0)=0,i=1,\cdots,N$.\par
For this purpose, we introduce some standard assumptions and Lemmas.\\
\textbf{Assumption 1.} The linear autonomous system is neutrally stable, that is, all the eigenvalues of $A$ are semi-simple with zero real parts.\\
\textbf{Assumption 2.} The digraph $\bar{\mathcal{G}}$ contains a directed spanning tree with node $v_0$ as its root.
\begin{lemma}\cite{Su-AC-2012}\label{lemma1}
Consider a weighted digraph $\bar{\mathcal{G}}=\{\bar{\mathcal{V}},\bar{\mathcal{E}},\bar{\mathcal{A}}\}$ with  $\bar{\mathcal{V}}=\{v_0,v_1,\cdots,v_N\}$, $\bar{\mathcal{E}}\subseteq \bar{\mathcal{V}}\times \bar{\mathcal{V}}$ and $\bar{\mathcal{A}}=[a_{ij}]\in\mathbb{R}^{(N+1)\times (N+1)}$. Let $\mathcal{L}=[l_{ij}]_{i,j=1}^N$ with $l_{ij}=-a_{ij},i\neq j$ and $l_{ii}=\sum_{j=1}^Na_{ij}, \Delta=\text{diag}\{a_{10},\cdots,a_{N0}\}$, and $\mathcal{H}=\mathcal{L}+\Delta$. Then, under Assumption 2, all the eigenvalues of $\mathcal{H}$ have positive real parts.
\end{lemma}
\begin{lemma}\cite{Tuna-arXiv-2008}\label{lemma3}
Let $A\in\mathbb{R}^{n\times n}$ and $C\in\mathbb{R}^{n\times m}$ satisfy
\begin{equation*}
AP+PA^T+I_n-PC^TCP=0
\end{equation*}
for some symmetric positive definite $P$. Then, for all $\sigma>1$ and $\omega\in\mathbb{R}$, matrix $A^T-(\sigma+j\omega)C^TCP$ is Hurwitz.
\end{lemma}

\begin{lemma}\label{lemma2}\cite{Hassan-NS-1996}
If $ \lim_{t\rightarrow\infty}\int_0^tf(\tau)d\tau$ exists and is finite, and $f(t)$ is a uniformly continuous function, then $\lim_{t\rightarrow\infty}f(t)=0$.
\end{lemma}
\textbf{Remark 1.}~The system $(\ref{systemF})$ is often called as the parametric strict feedback system in the literature. It is commonly encountered in many nonlinear control problems. Actually, under some mild conditions, a class of general nonlinear systems $\dot{\chi}=f_0(\chi)+\sum_{l=1}^p\theta_lf_l(\chi)+g(\chi)u,y=h(\chi)$ can be transformed into such form \cite{Kristic-NACD-1995}. Based on the adaptive backsteping approach, the adaptive output consensus problem of multi-agent systems in this form has been investigated in \cite{Wangwei-A-2014}-\cite{HuangJingshuai-A-2017}. Unfortunately, some restrictive conditions were imposed, e.g., each agent needs to know the state information and nonlinear functions of its neighbors \cite{Wangwei-A-2014}-\cite{Wangwei-A-2017}, or the filter information of its neighbors \cite{HuangJingshuai-A-2017}, and the system orders of all agents need to be the same. In our work, the system order needs not to be the same for all agents, and each agent only needs to know the output information of its neighbors. This considerably reduces communication burden and facilitates the implementation of the controller.

\section{Main Results}
\subsection{A novel distributed dynamic compensator}
In this section, a novel distributed dynamic compensator is developed to address the challenges caused by heterogeneous dynamics. First, we consider the following dynamic compensator:
\begin{eqnarray}
\dot{\eta}_{i,l}&=&A\eta_{i,l}-KC(\eta_{i,l}-\eta_{i,l+1}),l=1,\cdots,r_i,\nonumber\\
\dot{\eta}_{i,r_i+1}&=&A\eta_{i,r_i+1}-KC\sum\limits_{j=0}^Na_{ij}(\eta_{i,r_i+1}-\eta_{i,1})-Ke_{vi},\nonumber\\
\hat{y}_i&=&C\eta_{i,1}, i=1,\cdots,N,\label{compensator}
\end{eqnarray}
where $K\in\mathbb{R}^{v\times 1}$ is a constant matrix to be designed later. Let $\hat{e}_{i,1}=y_i-\hat{y}_i,i=1,\cdots,N$, and $\hat{e}=[\hat{e}_{1,1},\cdots,\hat{e}_{N,1}]^T$. Moreover, let $\bar{\eta}_{i,l}=\eta_{i,l}-x_0,i=1,\cdots,N,l=1,\cdots,r_i+1$. Then, for this dynamic compensator, we have the following results.
\begin{theorem}\label{theorem1}
Consider the dynamic compensator $(\ref{compensator})$ with $K$ being designed by $(\ref{gainK})$. Under Assumption 2, there exist a $\mathcal{KL}$-function $\beta$ and a $\mathcal{K}$-function $\gamma$ such that for $i=1,\cdots,N,l=1,\cdots,r_i+1,$
\begin{equation}
\|\bar{\eta}_{i,l}(t)\|\leq \beta(\bar{\eta}_{i,l}(0),t)+\gamma(\|\hat{e}\|_{\infty}), t\geq 0.
\end{equation}
In particular, if $\lim_{t\rightarrow \infty}\hat{e}_{i,1}=0,i=1,\cdots,N$, then $\lim_{t\rightarrow \infty}\bar{\eta}_{i,l}(t)=0,i=1,\cdots,N,l=1,\cdots,r_i+1.$
\end{theorem}
\textbf{Proof.}~From $(\ref{systemL})$ and $(\ref{compensator})$, we have
\begin{small}
\begin{eqnarray}
\dot{\bar{\eta}}_{i,l}&=&A\bar{\eta}_{i,l}-KC(\bar{\eta}_{i,l}-\bar{\eta}_{i,l+1}),l=1,\cdots,r_i,\nonumber\\
\dot{\bar{\eta}}_{i,r_i+1}\!&=&\!A\bar{\eta}_{i,r_i+1}\!-\!K\sum\limits_{j=0}^Na_{ij}(C\eta_{i,r_i+1}\!-\!C\eta_{i,1}+y_i-y_j).\label{eq1}
\end{eqnarray}
\end{small}
Observe that
\begin{eqnarray}
&&\sum\limits_{j=0}^Na_{ij}(C\eta_{i,r_i+1}-C\eta_{i,1}+y_i-y_j)\nonumber\\
&&~=\sum\limits_{j=0}^Na_{ij}(y_i-\hat{y}_i)+a_{i0}C\bar{\eta}_{i,r_i+1}\nonumber\\
&&~~~+\sum\limits_{j=1}^Na_{ij}(C\eta_{i,r_i+1}-C\eta_{j,1}+\hat{y}_j-y_j).\label{eq2}
\end{eqnarray}
Thus, submitting $(\ref{eq2})$ into $(\ref{eq1})$ yields
\begin{eqnarray}
&&\dot{\bar{\eta}}_{i,l}=A\bar{\eta}_{i,l}-KC(\bar{\eta}_{i,l}-\bar{\eta}_{i,l+1}),l=1,\cdots,r_i,\nonumber\\
&&\dot{\bar{\eta}}_{i,r_i+1}=A\bar{\eta}_{i,r_i+1}\!-\!a_{i0}KC\bar{\eta}_{i,r_i+1}
\!-\!\sum\limits_{j=1}^Na_{ij}KC(\bar{\eta}_{i,r_i+1}\!-\!\bar{\eta}_{j,1})\nonumber\\
&&~~~~~~~~-K\sum\limits_{j=0}^Na_{ij}\hat{e}_{i,1}+K\sum\limits_{j=1}^Na_{ij}\hat{e}_{j,1}.\label{eq3}
\end{eqnarray}
Let $\bar{\eta}_i=[\bar{\eta}_{i,1},\cdots,\bar{\eta}_{i,r_i+1}]^T$ and $\bar{\eta}=[\bar{\eta}_1^T,\cdots,\bar{\eta}_N^T]^T$. Moreover, we define the block matrix $\hat{\mathcal{L}}=[\hat{L}_{ij}]_{i,j=1}^N$ as
\begin{small}
\begin{eqnarray}
&&\hat{L}_{ii}=\left[\begin{array}{cccccc}
               1&-1&0&\cdots&0&0\\
               0&1&-1&\cdots&0&0\\
               \vdots&\vdots&\vdots&\ddots&\vdots&\vdots\\
               0&0&0&\cdots&1&-1\\
               0&0&0&\cdots&0&\sum\limits_{j=1}^Na_{ij}\end{array}\right]\in\mathbb{R}^{(r_i+1)\times (r_i+1)},\nonumber\\
&&\hat{L}_{ij}=\left[\begin{array}{ccc}
 0&\cdots&0\\
 \vdots&\ddots&\vdots\\
 -a_{ij}&\cdots&0
 \end{array}\right]\in\mathbb{R}^{(r_i+1)\times(r_j+1)},i\neq j.
\end{eqnarray}
\end{small}
Then, the system $(\ref{eq3})$ can be expressed as
\begin{equation}
\dot{\bar{\eta}}=[I_{(r_1+1)+\cdots+(r_N+1)}\otimes A-(\hat{\mathcal{L}}+\hat{\Delta})\otimes KC]\bar{\eta}+M\hat{e},\label{systemeta}
\end{equation}
where $M\in\mathbb{R}^{\nu(r_1+1+\cdots+r_N+1)\times N}$ is a constant matrix, and $\hat{\Delta}=\text{diag}\{\hat{\Delta}_1,\cdots,\hat{\Delta}_N\}$ with $\hat{\Delta}_i=\text{diag}\{0,\cdots,0,a_{i0}\}\in\mathbb{R}^{(r_i+1)\times (r_i+1)}$.\par
In what follows, we prove that there exists matrix $K$ such that $\hat{A}:=I_{(r_1+1)+\cdots+(r_N+1)}\otimes A-(\hat{\mathcal{L}}+\hat{\Delta})\otimes KC$ is Hurwitz. First, we define a weighted digraph $\hat{\bar{\mathcal{G}}}=(\hat{\bar{\mathcal{V}}},\hat{\bar{\mathcal{E}}},\hat{\bar{\mathcal{A}}})$ from $\bar{\mathcal{G}}$ according to the following rules:\\
$(1):$~The node set $\hat{\bar{\mathcal{V}}}$ is defined by replacing the vertexes $v_i,i=1,\cdots,N$ of $\bar{\mathcal{G}}$ with vertices $\hat{v}_{i}^{1},\cdots,\hat{v}_{i}^{r_i},\hat{v}_{i}^{r_i+1}$, and the vertex $v_0$ with vertex $\hat{v}_0$;\\
$(2):$~The edges $(\hat{v}_{i}^{r_i+1},\hat{v}_{i}^{r_i}),\cdots, (\hat{v}_{i}^{2},\hat{v}_{i}^{1}),i=1,\cdots,N$ contain in $ \hat{\bar{\mathcal{E}}}$;\\
$(3):$~The edges $(\hat{v}_{i}^{1},\hat{v}_{j}^{r_j+1})\in \hat{\bar{\mathcal{E}}},i,j=1,\cdots,N$ if and only if $(v_i,v_j)\in \mathcal{E}$, and the edges $(\hat{v}_{0},\hat{v}_{i}^{r_i+1})\in\hat{\bar{\mathcal{E}}},i=1,\cdots,N$ if and only $(v_0,v_i)\in \bar{\mathcal{E}}$;\\
$(4):$~The weighted adjacency
matrix $\hat{\bar{\mathcal{A}}}=[\hat{\mathcal{A}}_{ij}]_{i,j=0}^N$ takes the following forms:
\begin{eqnarray}
&&\hat{\mathcal{A}}_{ii}=\left[\begin{array}{cc}
                   0& I_{r_i}\\
                   0&0
                   \end{array}\right],i\neq 0, \hat{\mathcal{A}}_{00}=0\nonumber\\
&&\hat{\mathcal{A}}_{ij}=\left[\begin{array}{ccc}
 0&\cdots&0\\
 \vdots&\ddots&\vdots\\
 a_{ij}&\cdots&0
 \end{array}\right],i,j\neq 0,i\neq j,\nonumber\\
 && \hat{\mathcal{A}}_{0j}=\textbf{0}_{r_j+1}^T,\hat{\mathcal{A}}_{j0}=[0,\cdots,0,a_{0j}]^T,j=1,\cdots,N.\label{blockmatrix}
\end{eqnarray}\par
Therefore, the nodes $v_i,i=1,\cdots, N$ in $\bar{\mathcal{G}}$ is associated to the nodes $\hat{v}_i^1,\cdots,\hat{v}_i^{r_i+1},$ and the node $v_0$ is associated to the node $\hat{v}_0$ in $\hat{\bar{\mathcal{G}}}$. Since the digraph $\bar{\mathcal{G}}$ contains a directed spanning tree with
node $v_0$ as root, for each node $v_j,j=1,\cdots,N$, there exists a path from node $v_0$ to node $v_j$. Without loss of generality, we assume that  $(v_0,v_{i_1}),(v_{i_1},v_{i_2}),\cdots,(v_{i_k},v_j)$ is the directed path from node $v_0$ to node $v_j$. Then, from items $(2)$ and $(3)$ above, we conclude that  there exists a directed path $(\hat{v}_0,\hat{v}_{i_1}^{r_{i_1}+1})$,
$(\hat{v}_{i_1}^{r_{i_1}+1},\hat{v}_{i_1}^{r_{i_1}})$,
$\cdots,(\hat{v}_{i_1}^2,\hat{v}_{i_1}^1)$,
$(\hat{v}_{i_1}^1,\hat{v}_{i_2}^{r_{i_2}+1})$,
$\cdots,(\hat{v}_{i_k}^1,\hat{v}_{j}^{r_j+1})$,
$\cdots,(\hat{v}_{j}^2,\hat{v}_{j}^1)$ from node $\hat{v}_0$ to node $\hat{v}_j^{1}$ in $\hat{\bar{\mathcal{G}}}$. As a result, the digraph $\hat{\bar{\mathcal{G}}}$ contains a directed spanning tree with
node $\hat{v}_0$ as its root.
 Note that the matrix $\mathcal{\hat{L}}$ is the Laplacian matrix of the digraph $\hat{\mathcal{G}}=(\hat{\mathcal{V}},\hat{\mathcal{E}})$, where $\hat{\mathcal{V}}=\{\hat{v}_i^{l},i=1,\cdots,N,l=1,\cdots,r_i+1\}$, and $\hat{\mathcal{E}}\subseteq \hat{\mathcal{V}}\times \hat{\mathcal{V}}$ is obtained from $\hat{\bar{\mathcal{E}}}$ by removing all the edges between the node $\hat{v}_0$ and the nodes in $\hat{\mathcal{V}}$. Thus, from lemma \ref{lemma1}, all the eigenvalues of $\hat{\mathcal{H}}:=\hat{\mathcal{L}}+\hat{\Delta}$ have positive real parts. \par
Since $(A,C)$ is detectable, there exists a unique solution $P_0=P_0^T>0$ to the following Riccati equation:
\begin{equation}
AP_0+P_0A^T+I_v-P_0C^TCP_0=0.
\end{equation}
Let
\begin{equation}
K=\mu P_0C^T,\label{gainK}
\end{equation}
where $\mu\geq \frac{1}{\text{Re}\{\lambda_1(\hat{\mathcal{H}})\}}$, and $\lambda_1(\hat{\mathcal{H}})$ denotes the eigenvalue of $\hat{\mathcal{H}}$ with the smallest real part. Then, by Lemma \ref{lemma3}, $A^T-\lambda_i(\hat{\mathcal{H}})\mu C^TCP_0,i=1,\cdots,r_1+\cdots+r_N+N$, are Hurwitz. Thus, $\hat{A}$ is Hurwitz.\par
 Since $\hat{A}$ is Hurwitz, there exists a positive matrix $P_1=P_1^T$ such that $\hat{A}^TP_1+P_1\hat{A}=-I$. Let $V_{\bar{\eta}}=\bar{\eta}^TP_1\bar{\eta}$. Then, by Young's inequality $2ab\leq \frac{1}{2}a^2+2b^2$, the time derivative of $V_{\bar{\eta}}$ along the $\hat{\eta}$ system $(\ref{systemeta})$ satisfies
 \begin{eqnarray}
 \dot{V}_{\bar{\eta}}&=&-\bar{\eta}^T{\bar{\eta}}+2{\bar{\eta}}^TP_1M\hat{e}\nonumber\\
 &\leq&-\frac{1}{2}\|{\bar{\eta}}\|^2+2\|P_1\|^2\|M\|^2\|\hat{e}\|^2.
 \end{eqnarray}
 Thus, the $\bar{\eta}$ system is ISS with $\bar{\eta}$ as the state and $\hat{e}$ as input, $i.e,$ there exist  a $\mathcal{KL}$-function $\beta$ and a $\mathcal{K}$-function $\gamma$ such that for $i=1,\cdots,N,l=1,\cdots,r_i+1,$
\begin{equation}
\|\bar{\eta}_{i,l}(t)\|\leq \beta(\bar{\eta}_{i,l}(0),t)+\gamma(\|\hat{e}\|_{\infty}), t\geq 0.\label{eq4}
\end{equation}
In particular, if $\lim_{t\rightarrow \infty}\hat{e}=0$, from $(\ref{eq4})$, one have
\begin{equation}
\lim\limits_{t\rightarrow \infty}\bar{\eta}_{i,l}(t)=0,i=1,\cdots,N,l=1,\cdots,r_i+1.
\end{equation}
The proof is completed. $\hfill \blacksquare$\\
\textbf{Remark 2:} By means of the distributed dynamic compensator $({\ref{compensator}})$, the Problem 1 is converted into the global asymptotic tracking problem of nonlinear system $(\ref{systemF})$. Actually, we only need to design adaptive controllers $u_i,i=1,\cdots, N$ such that $\lim_{t\rightarrow \infty}(y_i-\hat{y}_i)=0,i=1,\cdots,N$. Then, from Theorem \ref{theorem1}, one have $\lim_{t\rightarrow \infty}(y_i-y_0)=\lim_{t\rightarrow \infty}(y_i-\hat{y}_i)+\lim_{t\rightarrow \infty}C(\eta_{i,1}-x_0)=0$.
\subsection{Adaptive controller design}
In this section, we design adaptive controllers for $(\ref{systemF})$ by means of adaptive backstepping approach. In what follows, to make the expression more
concise, let the error for each design step be defined as
\begin{eqnarray}
\hat{e}_{i,1}&=&y_i-\hat{y}_i,\nonumber\\
\hat{e}_{i,l}&=&x_{i,l}-\alpha_{i,l-1},l=2,\cdots,r_i,
\end{eqnarray}
where $\alpha_{i,l-1}$  is the virtual control signal to be designed at the step $l-1$. Let $\hat{\theta}_i$ be the estimate of $\theta_i$, and
 \begin{equation}
 \tilde{\theta}_i=\hat{\theta}_i-\theta_i.
 \end{equation}
 Besides, we shall employ positive scalars $c_{i,l},l=1,\cdots,r_i$ as design parameters in the subsequent design steps without restating.\par
 {\bf \textit{ Step $1$.}} From $(\ref{systemF})$ and $(\ref{compensator})$, the derivative of $\hat{e}_{i,1}$ is expressed as
 \begin{eqnarray}
 \dot{\hat{e}}_{i,1}=x_{i,2}+\psi_{i,1}(x_{i,1})^T\theta_i-CA\eta_{i,1}+CKC(\eta_{i,1}-\eta_{i,2}).\label{eq10}
 \end{eqnarray}
 Consider the Lyapunov function
 \begin{equation}
 V_{i,1}=\frac{1}{2}\hat{e}_{i,1}^2+\frac{1}{2}\tilde{\theta}_i^T\tilde{\theta}_i.
\end{equation}
 A direct calculation results in
\begin{eqnarray}
\dot{V}_{i,1}&=&\hat{e}_{i,1}\hat{e}_{i,2}+\hat{e}_{i,1}[\alpha_{i,1}+\psi_{i,1}^T\theta_i-CA\eta_{i,1}\nonumber\\
&&+CKC(\eta_{i,1}-\eta_{i,2})]+\tilde{\theta}_i^T\dot{\hat{\theta}}_i.
\end{eqnarray}
Choose the first tuning function $\tau_{i,1}$ and the virtual control signal $\alpha_{i,1}$ as
\begin{small}
\begin{eqnarray}
&&\tau_{i,1}=\psi_{i,1}\hat{e}_{i,1},\nonumber\\
&&\alpha_{i,1}\!=\!-\!c_{i,1}\hat{e}_{i,1}\!-\!\psi_{i,1}^T\hat{\theta}_{i}\!+\!CA\eta_{i,1}
\!-\!CKC(\eta_{i,1}\!-\!\eta_{i,2}).\label{controller-step1}
\end{eqnarray}
\end{small}
Then, it can be checked that
\begin{eqnarray}
\dot{V}_{i,1}=-c_{i,1}\hat{e}_{i,1}^2+\hat{e}_{i,1}\hat{e}_{i,2}+\tilde{\theta}_i^T(\dot{\hat{\theta}}_i-\tau_{i,1}).
\end{eqnarray}\par
 {\bf \textit{ Step $k(2\leq k\leq r_i-1)$.} }Similarly, the derivative of $\hat{e}_{i,k}$, by considering $(\ref{systemF})$ and $(\ref{compensator})$, can be computed as
 \begin{small}
\begin{eqnarray}
&&\dot{\hat{e}}_{i,k}=\hat{e}_{i,k+1}+\alpha_{i,k}+\psi_{i,k}^T\theta_i-\sum\limits_{l=1}^{k-1}\frac{\partial \alpha_{i,k-1}}{\partial x_{i,l}}(x_{i,l+1}+\psi_{i,l}^T\theta_i)\nonumber\\
&&~~\!-\!\sum\limits_{l=1}^{k}\!\frac{\partial \alpha_{i,k-1}}{\partial \eta_{i,l}^T}[A\eta_{i,l}\!-\!KC(\eta_{i,l}\!-\!\eta_{i,l+1})]
\!-\!\frac{\partial \alpha_{i,k-1}}{\partial \hat{\theta}_i^T}\dot{\hat{\theta}}_i.
\end{eqnarray}
\end{small}
Define the quadratic function $V_{i,k}$ as
\begin{equation}
V_{i,k}=V_{i,k-1}+\frac{1}{2}\hat{e}_{i,k}^2,
\end{equation}
where the derivative of $V_{i,k-1}$, by induction, satisfies
\begin{eqnarray}
&&\dot{V}_{i,k-1}= -\sum\limits_{l=1}^{k-1}c_{i,l}\hat{e}_{i,l}+\hat{e}_{i,k-1}\hat{e}_{i,k}+\tilde{\theta}_i^T(\dot{\hat{\theta}}_i-\tau_{i,k-1})\nonumber\\
&&~~-\sum\limits_{l=2}^{k-1}\hat{e}_{i,l}\frac{\partial \alpha_{i,l-1}}{\partial \hat{\theta}_i^T}(\dot{\hat{\theta}}_i-\tau_{i,k-1}).
\end{eqnarray}
A direct calculation  leads to
 \begin{eqnarray}
 &&\dot{V}_{i,k}=-\sum\limits_{l=1}^{k-1}c_{i,l}\hat{e}_{i,l}+\hat{e}_{i,k-1}\hat{e}_{i,k}+\tilde{\theta}_i^T(\dot{\hat{\theta}}_i-\tau_{i,k-1})\nonumber\\
&&~~-\sum\limits_{l=2}^{k-1}\hat{e}_{i,l}\frac{\partial \alpha_{i,l-1}}{\partial \hat{\theta}_i^T}(\dot{\hat{\theta}}_i-\tau_{i,k-1})+\hat{e}_{i,k}\hat{e}_{i,k+1}\nonumber\\
&&~~+\hat{e}_{i,k}\left[\alpha_{i,k}+\psi_{i,k}^T\theta_i-\sum\limits_{l=1}^{k-1}\frac{\partial \alpha_{i,k-1}}{\partial x_{i,l}}(x_{i,l+1}+\psi_{i,l}^T\theta_i)\right.\nonumber\\
&&~~\left.\!-\!\sum\limits_{l=1}^{k}\frac{\partial \alpha_{i,k-1}}{\partial \eta_{i,l}^T}\![A\eta_{i,l}\!-\!KC(\eta_{i,l}\!-\!\eta_{i,l+1})]
-\frac{\partial \alpha_{i,k-1}}{\partial \hat{\theta}_i^T}\dot{\hat{\theta}}_i\right].\nonumber
 \end{eqnarray}
 Then, choosing the tuning function $\tau_{i,k}$ and the virtual control signal $\alpha_{i,k}$ as
 \begin{small}
 \begin{eqnarray}
 &&\tau_{i,k}=\tau_{i,k-1}+\left(\psi_{i,k}-\sum\limits_{l=1}^{k-1}\psi_{i,l}\frac{\partial \alpha_{i,k-1}}{\partial x_{i,l}}\right)\hat{e}_{i,k},\nonumber\\
 &&\alpha_{i,k}\!=\!-c_{i,k}\hat{e}_{i,k}\!-\!\hat{e}_{i,k-1}\!-\!\psi_{i,k}^T\hat{\theta}_i+\sum\limits_{l=1}^{k-1}\frac{\partial \alpha_{i,k-1}}{\partial x_{i,l}}(x_{i,l+1}+\psi_{i,l}^T\hat{\theta}_i)\nonumber\\
 &&~~+\sum\limits_{l=1}^{k}\frac{\partial \alpha_{i,k-1}}{\partial \eta_{i,l}^T}[A\eta_{i,l}-KC(\eta_{i,l}-\eta_{i,l+1})]
 +\frac{\partial \alpha_{i,k-1}}{\partial \hat{\theta}_i^T}\tau_{i,k}\nonumber\\
 &&~~+\sum\limits_{l=2}^{k-1}\hat{e}_{i,l}\frac{\partial \alpha_{i,l-1}}{\partial \hat{\theta}_i^T}\left(\psi_{i,k}-\sum\limits_{l=1}^{k-1}\psi_{i,l}\frac{\partial \alpha_{i,k-1}}{\partial x_{i,l}}\right),
 \end{eqnarray}
 \end{small}
 one can obtain
 \begin{eqnarray}
 &&\dot{V}_{i,k}=-\sum\limits_{l=1}^{k}c_{i,l}\hat{e}_{i,l}^2+\hat{e}_{i,k}\hat{e}_{i,k+1}+\tilde{\theta}_i^T(\dot{\hat{\theta}}_i-\tau_{i,k})\nonumber\\
 &&~~~-\sum\limits_{l=2}^{k}\hat{e}_{i,l}\frac{\partial \alpha_{i,l-1}}{\partial \hat{\theta}_i^T}(\dot{\hat{\theta}}_i-\tau_{i,k}).\label{eq5}
 \end{eqnarray}\par
 {\bf \textit{ Step $r_i$.}}~The derivative of $\hat{e}_{i,r_i}$ is computed as
 \begin{small}
 \begin{eqnarray}
 &&\dot{\hat{e}}_{i,r_i}=u_i+\psi_{i,r_i}^T\theta_i-\sum\limits_{l=1}^{r_i-1}\frac{\partial \alpha_{i,r_i-1}}{\partial x_{i,l}}(x_{i,l+1}+\psi_{i,l}^T\theta_i)\nonumber\\
 &&~\!-\!\sum\limits_{l=1}^{r_i}\!\frac{\partial \!\alpha_{i,r_i-1}}{\partial \!\eta_{i,l}^T}\!(A\eta_{i,l}\!-\!KC(\eta_{i,l}\!-\!\eta_{i,l+1}))\!-\!\frac{\partial \alpha_{i,r_i-1}}{\partial \hat{\theta}_i^T}\dot{\hat{\theta}}_i.\label{eq6}
 \end{eqnarray}
 \end{small}
 The Lyapunov function is constructed as
 \begin{eqnarray}
 V_{i,r_i}=V_{i,r_i-1}+\frac{1}{2}\hat{e}_{i,r_i}^2,\label{LyapunovF}
 \end{eqnarray}
and the adaptive control law $u_i$ is chosen  as
\begin{eqnarray}
&&\dot{\hat{\theta}}_i=\tau_{i,r_n}:=\tau_{i,r_i-1}+\left(\psi_{i,r_i}-\sum\limits_{l=1}^{r_i-1}\psi_{i,l}\frac{\partial \alpha_{i,r_i-1}}{\partial x_{i,l}}\right)\hat{e}_{i,r_i},\nonumber\\
&&u_i=\alpha_{i,r_i}:=-c_{i,r_i}\hat{e}_{i,r_i}-\hat{e}_{i,r_i-1}-\psi_{i,r_i}^T\hat{\theta}_i\nonumber\\
&&~~~~+\sum\limits_{l=1}^{r_i-1}\frac{\partial \alpha_{i,r_i-1}}{\partial x_{i,l}}(x_{i,l+1}+\psi_{i,l}^T\hat{\theta}_i)+\frac{\partial \alpha_{i,r_i-1}}{\partial \hat{\theta}_i^T}\tau_{i,r_i}\nonumber\\
 &&~~~~+\sum\limits_{l=1}^{r_i}\frac{\partial \alpha_{i,r_i-1}}{\partial \eta_{i,l}^T}[A\eta_{i,l}-KC(\eta_{i,l}-\eta_{i,l+1})]
 \nonumber\\
 &&~~~~+\sum\limits_{l=2}^{r_i-1}\hat{e}_{i,l}\frac{\partial \alpha_{i,l-1}}{\partial \hat{\theta}_i^T}\left(\psi_{i,r_i}-\sum\limits_{l=1}^{r_i-1}\psi_{i,l}\frac{\partial \alpha_{i,r_i-1}}{\partial x_{i,l}}\right).\label{controller}
\end{eqnarray}
\textbf{Remark 3:} In \cite{Wangwei-A-2014}-\cite{Wangwei-A-2017}, each agent required constructing additional local estimates to account for the unknown parameters of its neighbors' dynamics. This inevitably results in a much complex controller. From $(\ref{controller})$, it is easy to know that our designed controller for each agent does not need  the additional local estimates.

 \subsection{Stability analysis}
 \begin{theorem}\label{theorem2}
 Consider the closed-loop system consisting of the $N$ nonlinear subsystems  $(\ref{systemF})$, the leader $(\ref{systemL})$, the dynamic compensator $(\ref{compensator})$ and the adaptive controllers $(\ref{controller})$. Under Assumptions 1 and 2, all signals in the closed-loop system are globally uniformly bounded, and asymptotic consensus tracking of all the subsystems' output to $y_0$ is achieved, $i.e.,\lim_{t\rightarrow \infty}[y_i(t)-y_0(t)]=0$ for $i=1,\cdots,N.$
 \end{theorem}
 \textbf{Proof.}~By using $(\ref{eq5})$ with $k=r_i-1$, we have
 \begin{eqnarray}
 &&\dot{V}_{i,r_i-1}=-\sum\limits_{l=1}^{r_i-1}c_{i,l}\hat{e}_{i,l}^2+\hat{e}_{i,r_i-1}\hat{e}_{i,r_i}+\tilde{\theta}_i^T(\dot{\hat{\theta}}_i-\tau_{i,r_i-1})\nonumber\\
 &&~~~~~-\sum\limits_{l=2}^{r_i-1}\hat{e}_{i,l}\frac{\partial \alpha_{i,l-1}}{\partial \hat{\theta}_i^T}(\dot{\hat{\theta}}_i-\tau_{i,r_i-1}).
 \end{eqnarray}
 Then, from $(\ref{eq6})$, the derivative of $V_{i,r_i}$ can be computed as
 \begin{eqnarray}
 &&\dot{V}_{i,r_i}=-\sum\limits_{l=1}^{r_i-1}c_{i,l}\hat{e}_{i,l}^2+\hat{e}_{i,r_i-1}\hat{e}_{i,r_i}+\tilde{\theta}_i^T(\dot{\hat{\theta}}_i-\tau_{i,r_i-1})\nonumber\\
 &&~~-\!\sum\limits_{l=2}^{r_i-1}\hat{e}_{i,l}\frac{\partial \alpha_{i,l-1}}{\partial \hat{\theta}_i^T}(\dot{\hat{\theta}}_i\!-\!\tau_{i,r_i-1})\!+\!\hat{e}_{i,r_i}\bigg[u_i\!-\!\alpha_{i,r_i}\!+\!\alpha_{i,r_i}\nonumber\\
 &&~~+\psi_{i,r_i}^T\theta_i-\sum\limits_{l=1}^{r_i-1}\frac{\partial \alpha_{i,r_i-1}}{\partial x_{i,l}}(x_{i,l+1}+\psi_{i,l}^T\theta_i)\!-\!\frac{\partial \alpha_{i,r_i-1}}{\partial \hat{\theta}_i^T}\dot{\hat{\theta}}_i\nonumber\\
 &&~~-\sum\limits_{l=1}^{r_i}\!\frac{\partial \!\alpha_{i,r_i-1}}{\partial \!\eta_{i,l}^T}\!(A\eta_{i,l}\!-\!KC(\eta_{i,l}\!-\!\eta_{i,l+1}))\bigg].\label{eq7}
 \end{eqnarray}
 Submitting $(\ref{controller})$ into $(\ref{eq7})$ results in
 \begin{eqnarray}
 &&\dot{V}_{i,r_i}=-\sum\limits_{l=1}^{r_i}c_{i,l}\hat{e}_{i,l}^2-\!\sum\limits_{l=2}^{r_i}\hat{e}_{i,l}\frac{\partial \alpha_{i,l-1}}{\partial \hat{\theta}_i^T}(\dot{\hat{\theta}}_i-\!\tau_{i,r_i})\!\nonumber\\
 &&~~~~+\tilde{\theta}_i^T(\dot{\hat{\theta}}_i-\tau_{i,r_i})+\!\hat{e}_{i,r_i}(u_i\!-\!\alpha_{i,r_i}).
 \end{eqnarray}
 Thus, according to the adaptive controller $(\ref{controller})$, one have
 \begin{eqnarray}
 \dot{V}_{i,r_i}=-\sum\limits_{l=1}^{r_i}c_{i,l}\hat{e}_{i,l}^2.\label{eq9}
 \end{eqnarray}
  Form $(\ref{LyapunovF})$, it is clear that $\hat{\theta}_i\in\mathcal{L}_{\infty}$ and $\hat{e}_{i,l}\in \mathcal{L}_{\infty}\cap\mathcal{L}_{2},l=1,\cdots,r_i.$ Then, from Theorem 1 and Assumption 1, one have $\bar{\eta}_{i,l}\in \mathcal{L}_{\infty}$ and $\eta_{i,l}\in \mathcal{L}_{\infty}$ for $i=1,\cdots,N,l=1,\cdots,r_i$. This implies $x_{i,l}\in \mathcal{L}_{\infty},i=1,\cdots,l=1,\cdots,r_i$. Therefore, all signals in the closed-loop system are globally uniformly bounded. Then, from $(\ref{eq10})$, we have $\dot{\hat{e}}_{i,1}\in \mathcal{L}_{\infty}$. In addition, according to Lemma \ref{lemma2} and the fact $\hat{e}_{i,1}\in \mathcal{L}_{\infty}\cap\mathcal{L}_{2},i=1,\cdots,N,$ one can verify that $\lim_{t\rightarrow \infty} \hat{e}_{i,1}=0,i=1,\cdots,N$. This together with Theorem 1, implies that $\lim_{t\rightarrow \infty}(y_i-y_0)=\lim_{t\rightarrow \infty}(y_i-\hat{y}_i)+\lim_{t\rightarrow \infty}C(\eta_{i,1}-x_0)=0$. The proof is completed. $\hfill \blacksquare$\\
 \textbf{Remark 4.} It should be noted that in \cite{Zhu-AC-2016}, the uncertain parameters must belong to a known compact set, and the controller for each subsystem needs to satisfy a class of small gain conditions. The small gain conditions may
result in sufficiently large control gains, and for some nonlinear systems with completely unknown parameters, it is unable or difficult to design controller satisfying small gain conditions. However, in our work, by means of the novel distributed dynamic compensator $(\ref{compensator})$, the controller $(\ref{controller})$ does not require the small gain conditions, and it hence can be adopted to the case that the uncertain parameters are completely unknown. \\
 \textbf{Remark 5.} By means of the distributed dynamic compensator $(\ref{compensator})$, we can use different control approaches to design tracking controller for each  subsystem. For example, for some agents with the following dynamic:
 \begin{eqnarray}
 &&\dot{x}_{i,l}=x_{i,l+1}+\psi_{i,l}(x_{i,1},\cdots,x_{i,l}),l=1,\cdots,r_i-1,\nonumber\\
&&\dot{x}_{i,r_i}=u_i+\psi_{i,r_i}(x_{i,1},\cdots,x_{i,r_i}),\nonumber\\
&&y_i=x_{i,1},\label{eqsystem}
 \end{eqnarray}
 where only the output $y_i$ can be measured by agent $i$,
 a linear-like output feedback controller can be designed for $(\ref{eqsystem})$ via the dynamic gain scaling technique \cite{ZhangandLin-A-2012}, which avoids the repeated derivatives of the nonlinearities depending on the observer states and the dynamic gain in backstepping approach.
 Thus, our proposed methodology can be applied to the leader-following output consensus problem of heterogeneous nonlinear multi-agent systems with non-identical structure. In particular, it can be used to deal with
  the output consensus problem of heterogeneous nonlinear multi-agent systems with unknown parameters and unknown non-identical control directions. Actually, by combination of adaptive backstepping technique and Nussbaum-type function, one can design an adaptive controller for each subsystem such that $\lim_{t\rightarrow \infty}(y_i-\hat{y}_i)=0$. Specifically, for the following heterogeneous nonlinear multi-agent systems with unknown and non-identical control directions:
\begin{eqnarray}
&&\dot{x}_{i,l}=x_{i,l+1}+\psi_{i,l}(x_{i,1},\cdots,x_{i,l})^T\theta_{i},l=1,\cdots,r_i-1,\nonumber\\
&&\dot{x}_{i,r_i}=b_iu_i+\psi_{i,r_i}(x_{i,1},\cdots,x_{i,r_i})^T\theta_{i},\nonumber\\
&&y_i=x_{i,1},i=1,\cdots,N,
\end{eqnarray}
where $b_i$ is a nonzero constant with unknown sign,  the adaptive controller for each subsystem is designed as
\begin{eqnarray}
&&u_i=-N_i(k_i)\alpha_{i,r_i},\nonumber\\
&&\dot{k}_i=-\hat{e}_{i,r_i}\alpha_{i,r_i},
\end{eqnarray}
where $N_i(k_i)$ is a Nussbaum function and $\alpha_{i,r_i}$ is defined in $(\ref{controller})$. Following the  the similar
proof in Theorem \ref{theorem2} with a minor modification, one can prove that all signals in the closed-loop system are globally bounded, and $\lim_{t\rightarrow \infty}(y_i-y_0)=0.$

\section{An illustrative example}
In this section, we consider a heterogeneous nonlinear multi-agent system connected by a communication graph shown in Figure 1, where the weighted adjacency matrix $\bar{\mathcal{A}}$ satisfies $a_{ij}=1$ if and only if $(v_j,v_i)\in\bar{\mathcal{E}}$.
\begin{figure}\centering
\includegraphics[width=0.16\textwidth]{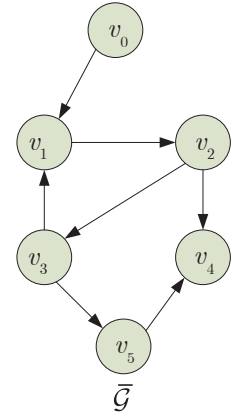}
\caption{The communication graph.} \label{fig1}
\end{figure}
The system is composed of agents with unknown parameters. In particular, agents $i,i=1,2,3$ are described by
\begin{eqnarray}
\dot{x}_{i,1}&=&x_{i,2}+x_{i,1}^2\theta_i,\nonumber\\
\dot{x}_{i,2}&=&u_i+\sin(x_{i,2})\theta_i,\nonumber\\
y_i&=&x_{i,1},
\end{eqnarray}
while agents $i,i=3,5$ are described by
\begin{eqnarray}
\dot{x}_{i,1}&=&u_i+\cos(x_{i,1})\theta_i,\nonumber\\
y_i&=&x_{i,1}.
\end{eqnarray}
The leader's signal $y_0$ is generated by the linear system $(\ref{systemL})$ with
\begin{eqnarray}
A=\left[\begin{array}{cc}
0&1\\
-1&0\end{array}\right],C=\left[\begin{array}{cc}
1&0
\end{array}\right].
\end{eqnarray}\par
Obviously, it can be seen that $(A,C)$ is detectable, and Assumptions 1-2 are satisfied. Therefore, by Theorem \ref{theorem2}, we can design a distributed adaptive controller of form $(\ref{control})$ for all agents such that all signals in the closed-loop system are globally uniformly bounded and $\lim_{t\rightarrow \infty}(y_i-y_0)=0,i=1,\cdots,N$. Following the
design procedure in Section III,  one can design the distributed adaptive controller as
follows. \par
\begin{itemize}
\item \emph{Step 1.} The distributed dynamic compensator is given in $(\ref{compensator})$ with $N=5,r_1=r_2=r_3=2,r_4=r_5=1,$ and
    \begin{equation}
    K=\left[\begin{array}{cc}
    17.3081 &5.3019
    \end{array}\right]^T.
    \end{equation}
\item \emph{Step 2.} For agents $i=1,2,3$, the first error $\hat{e}_{i,1}=x_{i,1}-\hat{y}_i$, the virtual control signal $\alpha_{i,1}$ and the first tuning function $\tau_{i,1}$ are given by $(\ref{controller-step1})$. The second error $\hat{e}_{i,2}$, the update law $\dot{\hat{\theta}}_i$ and the controller law $u_i$ are given, respectively, by
    \begin{small}
    \begin{eqnarray*}
    &&\hat{e}_{i,2}=x_{i,2}-\alpha_{i,1},\nonumber\\
    &&\dot{\hat{\theta}}_i=\tau_{i,1}+\left(\sin(x_{i,2})-x_{i,1}^2\frac{\partial \alpha_{i,1}}{\partial x_{i,1}}\right)\hat{e}_{i,2},\nonumber\\
    &&u_i=-c_{i,2}\hat{e}_{i,2}-\hat{e}_{i,1}-\sin(x_{i,2})\hat{\theta}_i+\frac{\partial \alpha_{i,1}}{\partial x_{i,1}}(x_{i,2}+x_{i,1}^2\hat{\theta}_i)\nonumber\\
    &&~~~~+\frac{\partial \alpha_{i,1}}{\partial \eta_{i,1}^T}[A\eta_{i,1}-KC(\eta_{i,1}-\eta_{i,2})]+\frac{\partial \alpha_{i,1}}{\partial \hat{\theta}_i}\tau_{i,2}\nonumber\\
    &&~~~~+\frac{\partial \alpha_{i,1}}{\partial \eta_{i,2}^T}[A\eta_{i,2}-KC(\eta_{i,2}-\eta_{i,3})].
   \end{eqnarray*}
   \end{small}

\item \emph{Step 3.} For agents $i=4,5$, the error $\hat{e}_{i,1}$, the update law $\dot{\hat{\theta}}_i$, and the controller law $u_i$ are given, respectively, by
    \begin{small}
        \begin{eqnarray*}
        &&\hat{e}_{i,1}=x_{i,1}-\hat{y}_{i},\nonumber\\
        &&\dot{\hat{\theta}}_i=\cos(x_{i,1})\hat{e}_{i,1},\nonumber\\
        &&u_i
    \!=\!-\!c_{i,1}\hat{e}_{i,1}\!-\!\cos(x_{i,1})\hat{\theta}_i\!+\!CA\eta_{i,1}\!-\!CKC(\eta_{i,1}\!-\!\eta_{i,2}).
        \end{eqnarray*}
        \end{small}
\end{itemize}\par
Simulation is performed with $c_{i,1}=c_{i,2}=c_{j,1}=1,i=1,2,3,j=4,5$,  $\theta_1=2.5,\theta_2=1.2,\theta_3=-2,\theta_4=-1,\theta_5=0.5,$ and the following initial conditions:
\begin{small}
\begin{eqnarray*}
&&x_1(0)=[0.1,-0.2]^T,x_2(0)=[0.5,1.2]^T,x_3(0)=[-2,1]^T,\nonumber\\
&&x_4(0)=-0.5,x_5(0)=0.25,x_0(0)=[1,-1],\hat{\theta}_1(0)=1.2,\nonumber\\
&&\hat{\theta}_2(0)=-1,\hat{\theta}_3(0)=0.5,\hat{\theta}_4(0)=0.2,\hat{\theta}_5(0)=-0.75,\nonumber\\
&&\eta_{1,1}(0)\!=\![0.1,0.2]^T,\eta_{1,2}(0)\!=\![1,-1.5]^T,\eta_{1,3}(0)\!=\![-1,-0.2]^T,\nonumber\\
&&\eta_{2,1}(0)\!=\![1\!,\!-0.5]^T,\eta_{2,2}(0)\!=\![-0.25,0.3]^T,\eta_{2,3}(0)\!=\![0.5,0.2]^T,\nonumber\\
&&\eta_{3,1}(0)\!=\![0.5,-0.4]^T,\eta_{3,2}(0)\!=\![0.6,-1]^T,\eta_{3,3}(0)\!=\![3,-0.2]^T,\\
&&\eta_{4,1}(0)=[2,-1.4]^T,\eta_{4,2}(0)=[2,1]^T,\\
&&\eta_{5,1}(0)=[1,2]^T,\eta_{5,2}(0)=[0.5,-0.75]^T.
\end{eqnarray*}
\end{small}
The simulation results are shown in Figure 2, which shows the effectiveness of  the design
methodology.
\begin{figure}\centering\label{fig2}
\includegraphics[width=0.42\textwidth]{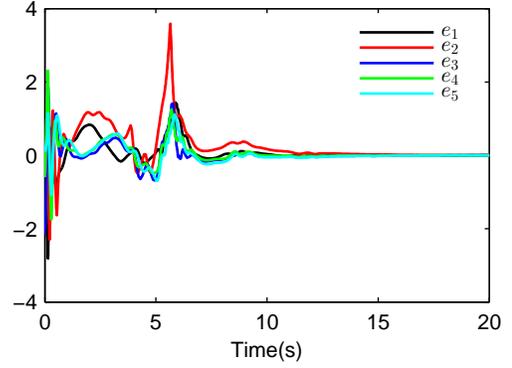}
\caption{Tracking errors of the agents.}
\end{figure}

\section{Conclusion}
The adaptive output consensus problem has been investigated in this note for a class of heterogeneous nonlinear multi-agent systems with unknown parameters. A novel distributed dynamic compensator has been developed to address the challenges caused by heterogeneous dynamics. The distributed dynamic compensator only requires the output information to be exchanged through communication networks. In addition, it can convert the original adaptive consensus problem into the problem of global asymptotic tracking for a class of nonlinear systems with unknown parameters. By means of adaptive backstepping approach, we have developed an adaptive tracking controller for each subsystem, which does not require the small gain conditions as in \cite{Zhu-AC-2016}. It has been proved that all signals in the closed-loop system are globally uniformly bounded, and the proposed scheme enables the outputs of all subsystems to track the output of leader asymptotically.

\end{document}